\documentstyle[]{article}
\newtheorem{thm}{Theorem}[section]

\newtheorem{cor}[thm]{Corollary}
\newtheorem{prop}[thm]{Proposition}
\newtheorem{conj}[thm]{Conjecture}
\newtheorem{dfn}[thm]{Definition}
\newtheorem{ex}[thm]{Example}

\newtheorem{rem}[thm]{Remark}

\newcommand{\F}{{\cal F}}
\newcommand{\del}{\partial}

\newcommand{\ra}{\rightarrow}

\font\bbb=msbm10 
\newcommand{\real}{\mbox{\bbb R}}

\newcommand{\rats}{\mbox{\bbb Q}}       
\newcommand{\zed}{\mbox{\bbb Z}}        
\newcommand{\inv}{^{-1}}

\newcommand{\curl}{\nabla\times}

\newcommand{\norm}[1]{\|#1\|}

\newcommand{\rest}[2]{\left. #1\right\vert_{#2}}                
\newcommand{\pf}{{\em Proof: }}

\newcommand{\qed}{$\Box$}

\newcommand{\eg}{{\em e.g.}}
\newcommand{\ie}{{\em i.e.}}
\newcommand{\cf}{{\em cf. }}
\begin{document}
%
%
\begin{center}
\Large 
{\bf Contact Topology and Hydrodynamics}
\normalsize
\vspace{0.1in}

%
%
John Etnyre\footnote{ E-mail : {\tt etnyre@math.utexas.edu}}
 and Robert Ghrist\footnote{ E-mail : {\tt ghrist@math.utexas.edu}} \\ 
Department of Mathematics \\
University of Texas, Austin \\
Austin, TX  78712 \\
\vspace{0.5in}

\end{center}
%
%
\begin{abstract}
We draw connections between the field of contact topology 
and the study of Beltrami fields in hydrodynamics
on Riemannian manifolds in dimension three. We demonstrate an 
equivalence between Reeb fields (vector fields which preserve a 
transverse nowhere-integrable plane field) up to scaling 
and rotational Beltrami fields on three-manifolds. Thus, we
characterise Beltrami fields in a metric-independant manner.

This correspondence yields a hydrodynamical reformulation 
of the Weinstein Conjecture, 
whose recent solution by Hofer (in several cases) implies 
the existence of closed orbits for all $C^\infty$ rotational 
Beltrami flows on $S^3$.
This is the key step for a positive solution to the hydrodynamical 
Seifert Conjecture: all $C^\omega$ steady state 
flows of a perfect incompressible fluid on $S^3$ possess closed
flowlines. In the case of Euler flows on 
$T^3$, we give general conditions for closed flowlines 
derived from the homotopy data of the normal bundle to the flow. 

\vspace{0.1in}
%
%
\end{abstract}

\section{Introduction}

The goal of this paper is to establish close relationships 
between the fields of hydrodynamics on Riemannian manifolds 
and contact topology in odd dimensions. For convenience and
the sake of applications, we restrict to dimension three, although
the basic relationships remain true in more general settings. 

In hydrodynamics, we consider the 
class of vector fields whose behavior is most fascinating and whose 
analysis has been most incomplete: these are the {\em Beltrami 
fields}, or fields which are parallel to their own curl. All such 
fields are solutions to Euler's equations of motion for a perfect 
incompressible fluid. Flows generated by such fields have 
several noteworthy properties, such as extremization of an 
energy functional (Equation~\ref{eq_Energy}), 
as well as the potential to display the phenomenon of  
{\em Lagrangian turbulence}, in which a volume-preserving flow   
has flowlines which fill up regions of space ergodically \cite{Ar2}. 

We relate the study of such flows with the rapidly developing 
field of contact topology. On a three-manifold, a contact structure 
is a maximally nonintegrable plane field. Though long a fixture in the 
literature on geometric classical mechanics (dating back at least
to the work of Lie),
it is only recently that topologists have made great progress in  
classifying such objects. The present state of affairs, brought
about via the work of Eliashberg, Gromov, Hofer, and others 
\cite{El1,El2,El3,Gr,Ho},
encompasses several very strong results in this area.
In particular, the method of analysing contact structures via 
{\em Reeb fields}, or transverse vector fields whose flow
preserve the contact form, has recently proven useful 
\cite{Ho,HWZ}.

In \S\ref{sec_Hydro} and \S\ref{sec_Contact}, 
we provide the necessary background 
information from each field, since we hope to promote
interaction between both areas of mathematics. 
Then, in \S\ref{sec_R&B}, 
we describe the correspondence between Beltrami fields and 
Reeb fields: 

{\em
{\bf  Theorem} The class of (nonsingular) vector fields on a 
three-manifold parallel to its (nonsingular) curl 
is identical to the class of Reeb fields under rescaling.
}

This yields a reformulation of the Weinstein Conjecture from 
symplectic topology into a hydrodynamical context: namely, 
that Beltrami flows on closed Riemannian three-manifolds 
must have closed flowlines. We spell out this reformulation 
in \S\ref{sec_WC}. 

The question of whether vector fields on three-manifolds are
forced to exhibit periodic solutions has a rich history
(see, \eg, \cite{Ku1,Ho}). The following 
conjecture of Seifert has been a focus of inquiry:

{\em
{\bf The Seifert Conjecture} Every $C^k$ vector field on 
$S^3$ possesses a closed orbit.
}

This conjecture was shown to be false for $k=1$ by 
Schweitzer \cite{Sc}, for $k=2$ by Harrison \cite{Ha}, and
for $k=\omega$ (\ie, real-analytic) by K. Kuperberg \cite{Ku1}. 
A recent theorem of Hofer \cite{Ho} resolves the Weinstein 
Conjecture on $S^3$. Via the correspondence between Beltrami
and Reeb fields, we derive the startling corollary that  
all $C^\infty$ rotational Beltrami flows on the three-sphere 
(and certain other manifolds) possess closed flowlines, no 
matter what the Riemannian metric implicated. 
Using this corollary, we provide a positive solution to the
Seifert Conjecture for perfect incompressible fluids:

{\em
{\bf Theorem} Every $C^\omega$ steady-state flow 
of a perfect incompressible fluid on $S^3$ possesses closed 
flowlines.
}

In other words, the $C^\omega$ Kuperberg plug cannot be 
parallel to its curl under {\em any} metric. Athough for 
Beltrami flows, we can relax the smoothness condition to 
$C^\infty$, it appears very difficult to improve the 
above theorem to $C^\infty$. For the non-Beltrami flows, 
we employ techniqes from singularity theory which are
dependant upon analyticity.


Finally, in \S\ref{sec_ABC}, we apply contact-topological methods
to a fundamental class of examples in hydrodynamics: 
spatially periodic flows on $\real^3$ (\ie, flows on $T^3$). 
In the case of Euclidean geometry, these include the 
{\em ABC flows}. By using a classification theorem of
Giroux \cite{Gi}, we give conditions which force the existence
of closed contractible orbits on $T^3$ in terms of the algebraic topology
of the vector field:

{\bf Theorem} {\em 
Any $C^\infty$ steady-state rotational Beltrami field on $T^3$ which is
homotopically nontrivial must have a contractible closed flowline. 
More generally, any $C^\omega$ steady-state Euler flow on $T^3$ 
which is homotopically nontrivial must have a closed flowline.
}

It is conjectured that the above theorem holds for {\em all} $C^\infty$ 
rotational Beltrami fields on $T^3$. It certainly does not 
hold for all Eulerian fluid flows, as irrational linear flow on 
$T^3$ is Eulerian (though homotopicaly trivial).

Throughout this paper, we use the language of differential 
forms, as it provides a convenient common basis for 
dealing with both hydrodynamics and contact geometry. 
Specifically, we use $\L_X$ to denote the Lie derivative 
along the vector field $X$, and we use $\iota_X$ to denote
contraction by $X$. 

\subsection{Topological hydrodynamics}\label{sec_Hydro}
The principal business of 
hydrodynamics is to understand the dynamical 
properties of fluid flows, the simplest class of which are
incompressible, inviscid flows described by the Euler equation.
The following treatment of hydrodynamics on Riemannian 
manifolds is based on the approach of Arnold and 
Khesin \cite{AK}, in which the reader will find 
several excellent references. 
In order to apply our results to the 
widest possible class, we adopt the following convention
(\cf \cite{AK}).
\begin{dfn}\label{def_Curl} {\em
The {\em curl} of a vector field $X$ on a Riemannian 3-manifold
$M$ with metric $g$ and a distinguished volume form $\mu$ is 
the unique vector field $\curl X$ given by
\begin{equation}
\iota_{(\curl X)}\mu = d\iota_Xg .
\end{equation}
}\end{dfn}
In the fluids literature, the curl of a velocity field is the 
{\em vorticity field}. 
Note that in the case where $\mu$ is the volume form for the 
metric $g$, the definition of $\curl X$ assumes the more familiar 
form \cite{AMR}
\begin{equation}
 \curl X = \left(*d\iota_Xg\right)^\# ,
\end{equation}
where $*$ denotes the Hodge star and $^\#$ denotes the 
isomorphism from 1-forms to vector fields derived from $g$.
The uniqueness of $\curl X$ in the above definition comes from 
the fact that for a fixed volume form $\mu$, the map 
$\iota_\cdot\mu$ is an isomorphism from vector fields to 2-forms. 

\begin{dfn}\label{def_Euler} {\em
Let $X(t)$ denote a (time-dependant) vector field on a Riemannian 
$3$-manifold $M$ with metric $g$ and distinguished volume
form $\mu$. Then $X(t)$ satisfies {\em Euler's equation for a 
perfect incompressible fluid} if 
\begin{equation}\label{eq_Euler}
\begin{array}{rcl}
\frac{\displaystyle \del X}{\displaystyle \del t} - X\times(\curl X) &=& 
        \nabla p \\
L_X\mu &=& 0 ,
\end{array}
\end{equation}
for some function $p(t)$ (the reduced pressure). 
The quantity $\L_X\mu$ vanishes if and only if the flow 
associated to $X$ is volume-preserving with respect to $\mu$.
We call vector fields (or corresponding flows) satisfying
(\ref{eq_Euler}) 
{\em Euler fields} (or {\em Euler flows} respectively).
}\end{dfn}

The field of hydrodynamics has been most successful in 
understanding steady Euler flows, or flows without time 
dependance. {\em For the remainder of this work, all vector fields
will be assumed to be steady-state.} The topology
of a steady state Euler flow is almost always very simple:
\begin{thm}[Arnold \cite{Ar1}]\label{thm_Arnold}
Let $X$ be a $C^\omega$ nonsingular Euler field on a closed Riemannian
three-manifold $M$. Then, if $X$ is not everywhere parallel to its curl, 
there exists a compact analytic subset $\Sigma\subset M$ of codimension 
at least one which splits $M$ into a finite collection of cells 
$T^2\times\real$. Each $T^2\times\{x\}$ is an invariant set for $X$
having flow conjugate to linear flow.  
\end{thm}
The proof is straightforward. The reduced 
pressure term $p$ is an integral for the flow since
$X\cdot\nabla p = -X\cdot(X\times(\curl X)) = 0$ (this is 
Bernoulli's Theorem).
This integral is nonzero 
if and only if the vorticity and velocity fields are independent.
In the nondegenerate case, for regular values
of $p$, the preimage is an invariant compact two-manifold
possessing a nonsingular flow: $\coprod T^2$. The preimage of the 
(finite number of) critical values of $p$ forms the singular set
$\Sigma$. By analyticity, $\Sigma$ may not contain open sets; 
hence, the topology of the flowlines are almost everywhere 
highly constrained {\em except} in the case where the curl of 
$X$ is colinear with $X$. 

The class of fields for which the above integral degenerates 
is extremely important.
\begin{dfn}\label{def_Beltrami} {\em
A vector field $X$ is said to be a {\em Beltrami field} if 
it is parallel to its curl: \ie, $\curl X = f X$ for some function 
$f$ on $M$. A {\em rotational} Beltrami field is one for which
$f\neq 0$; that is, it has nonzero curl.
}\end{dfn}

Beltrami fields have generated significant interest in the fluids 
community \cite{Ar2,D+,MP,Mo,FGV}. 
As shown by Arnold \cite{Ar2}, Beltrami fields 
extremize an energy functional within the class of flows 
conjugate by $\mu$-preserving diffeomorphisms. Furthermore, 
the flowlines of a Beltrami flow on a three-manifold are not always 
constrained to lie on a 2-torus, as is the typical case for non-Beltrami 
steady state Euler flows. Hence, the manifestation of 
``Lagrangian turbulence'' can only appear within the class of 
Beltrami fields. 

The classic examples of Beltrami fields which appear to exhibit 
Lagrangian turbulence are the {\em ABC flows}, generated by a certain 
family of vector fields on $T^3$ which are eigenfields of the curl operator
(see Equation~\ref{eq_ABC} in \S\ref{sec_ABC}). 
Although the ABC flows have been repeatedly analysed 
\cite{Ar2,D+,MP,FGV}, few results are known, 
other than for a handful of near-integrable examples. 
Beltrami fields in general are even less well understood. 
Typically, one wants to restrict to, say, Beltrami fields 
on Euclidean $\real^3$ or $T^3$, under the fixed standard 
metric and volume form. Under these restrictions, it becomes 
nearly impossible to do any sort of analysis on the class 
of Beltrami fields. A small perturbation, even within the class 
of $\mu$-preserving fields, almost always destroys the Beltrami 
property. We note in particular the difficulty of answering 
global questions about Beltrami flows, such as the existence 
of closed orbits, the presence of hydrodynamic instability, and the 
minimization of the energy functional. In this work, we
provide some topological tools which may prove robust enough
to overcome these difficulties.

\subsection{Contact topology}\label{sec_Contact}

For the sake of concreteness and applicatibility, we will 
restrict all definitions and discussions to the case of contact 
structures on three-manifolds, noting that several features hold 
on arbitrary odd-dimensional manifolds. For introductory 
treatments, see \cite{MS,Ae}. 

A {\em contact structure} on a three-manifold $M$ is a 
maximally nonintegrable plane field. That is, to each point $p\in M$,
we assign a hyperplane in $T_pM$, varying smoothly with $p$ in 
such a manner that the Frobenius condition fails everywhere. 
In particular, a contact structure is locally twisted at every point 
and may be thought of as an ``anti-foliation'' --- no disc may 
be embedded whose tangent planes agree with the plane field. 
\begin{dfn}{\em \label{def_Contact}
A {\em contact form} on $M^3$ is a one-form $\alpha$ on $M$ 
such that $\alpha\wedge d\alpha\neq 0$. That is, 
$\alpha\wedge d\alpha$ defines a 
volume form on $M$. A {\em contact structure} is a plane field
which is the kernel of a (locally defined) contact form:
\begin{equation}
\xi = {\mbox{ker}}(\alpha) = \{ v\in T_pM : \alpha(v)=0, p\in M \}.
\end{equation}
}\end{dfn}
It is usually sufficient to consider contact structures which are
the kernel of a globally defined contact 1-form: these are called 
{\em cooriented} contact structures. Of course, a double-cover suffices 
to make any plane field coorientable.

Contact structures have often been described as ``symplectic
structures in odd dimensions'' \cite{Ar3}, due to the fact that 
the restriction of $d\alpha$ to the plane field $\xi$ is a closed 
nondegenerate 2-form. However, given the definition in 
terms of nonintegrable plane fields, it is easy to see several 
elementary properties of contact structures. For example, 
unlike foliations, contact structures are structurally stable,
in the sense that perturbing $\alpha$ cannot force $\xi$ to 
be integrable. On the other hand, the manifestation of the Darboux 
Theorem in this context implies that every contact structure
locally looks like (is contactomorphic to, or diffeomorphic 
via a map which carries the contact structure to) the kernel
of $dz+x dy$ on $\real^3$ (see \cite{MS}).
Note the similarity with codimension-one foliations, which are 
locally equivalent to the kernel of $dz$ on $\real^3$. 

Again, as in foliation theory, the global features of a contact
structure are closely related to those of the manifold in which 
it sits. The classification of contact structures follows 
along lines similar to the Reeb-component versus taut 
perspective in (codimension-one) foliation theory \cite{ET}.

\begin{dfn}{\em 
Given a three-manifold $M$ with contact structure $\xi$, 
let $F\subset M$ be an embedded surface. Then the 
{\em characteristic foliation} on $F$, $F_\xi$, is the 
foliation on $F$ generated by the (singular) line field
\[
        {\cal F} = \left\{ T_pF \cap \xi_p : p\in F \right\} .
\]
A contact structure $\xi$ is {\em overtwisted} if there exists an 
embedded disc $D\subset M$ such that the characteristic 
foliation $D_\xi$ has a limit cycle. A contact structure 
which is not overtwisted is called {\em tight}. 
}\end{dfn}

The classification of overtwisted structures up to contactomorphism 
coincides with the classification of plane fields up to homotopy 
\cite{El1} and hence reduces to a problem in algebraic topology. 
The classification of tight structures, on the other hand, is 
far from complete. Like taut foliations, tight contact structures
exhibit several ``rigid'' features which make them relatively 
rare. The classification of tight contact structures is completed 
only on $S^3$ \cite{El3}, on certain lens spaces \cite{Et2}, and on 
the three-torus $T^3$ \cite{Gi}. It is unknown 
whether every closed orientable 3-manifold has a tight 
contact structure on it.

Given a contact form $\alpha$ generating the contact structure 
$\xi$, we may associate to it a vector field which is transverse
to $\xi$ and preserves $\alpha$ under the induced flow. 
Such vector fields were first considered by Reeb \cite{Re}.
\begin{dfn}{\em \label{def_Reeb}
Given a contact form $\alpha$ on $M$, the {\em Reeb field} 
associated to $\alpha$ is the unique vector field $X$ such that
\begin{equation}
        \iota_Xd\alpha = 0 \hspace{0.2in}
        {\mbox{  and  }} \hspace{0.2in}
        \iota_X\alpha =1 .
\end{equation}
}\end{dfn}
The condition that $\iota_X\alpha=1$ is a normalisation condition 
which corresponds to a time reparametrisation of the flow. As
we are primarily concerned with the topology of the flowlines, 
which does not depend on the parametrisation, we will also 
consider the class of {\em Reeb-like} fields, for which 
$\iota_Xd\alpha=0$ and $\iota_X\alpha>0$. 

The relationship between the dynamics of a Reeb field and 
the topology of the transverse contact structure have been 
analysed most notably by Hofer {\em et al.} 
\cite{Ho,HWZ}. We will consider these results in detail in 
\S\ref{sec_WC}.

\begin{ex}\label{ex_TightS3}{\em 
The standard contact structure on the unit $S^3\subset\real^4$ 
is given by the kernel of the 1-form 
\begin{equation}
        \alpha = \frac{1}{2}\left(
        x_1dx_2 - x_2dx_1 + x_3dx_4 - x_4dx_3 \right) .
\end{equation}
The Reeb field associated to $\alpha$ is the unit tangent field
to the standard Hopf fibration of $S^3$; hence, the contact 
structure $\xi$ can be visualised as the plane field orthogonal 
to the Hopf fibration (orthogonal with respect to the metric 
on the unit 3-sphere induced by the standard metric on $\real^4$).
It is a foundational result of Bennequin \cite{Be} that this
structure is tight; furthermore, by Eliashberg \cite{El3}, 
this is the unique tight contact structure on $S^3$ up to 
contactomorphism.
}\end{ex}

\section{Geometric properties of Beltrami flows}
\label{sec_R&B}

It is not coincidental that the Beltrami condition on a vector 
field --- that the flow must continually twist about itself --- is
reminiscent of the notion of a nonintegrable, everywhere
twisting plane field. We obtain a general equivalence 
between Beltrami and Reeb fields for arbitrary three-manifolds 
by working with moving frames and metrics adapted to the flow. 

\begin{thm}\label{thm_BelReeb}
Let $M$ be a Riemannian three-manifold. Any rotational 
Beltrami field on $M$ is a Reeb-like field
for some contact form on $M$. Conversely, given a contact form 
$\alpha$ on $M$ with Reeb field $X$, any nonzero rescaling of $X$
is a rotational Beltrami field for some Riemannian metric on $M$.
\end{thm}
\pf
Assume that $X$ is a Beltrami field where $\curl X=f X$ for
some $f>0$. Let $g$ denote the metric and $\mu$ the volume
form on $M$. On charts for $M$, choose an orthonormal 
frame $\{e_i\}_1^3$ such that $e_1=X/\norm{X}$. 
Note that we must work on charts only in the case where the 
Euler class of the vector field ($e(X)\in H^2(M;\zed)$) is 
nonzero. Denoting by 
$\{e^i\}_1^3$ the dual 1-form basis, we have that 
$\iota_Xg=\norm{X}e^1$. Let $\alpha$ denote the 
one-form $\iota_Xg=\norm{X}e^1$, which is globally defined since 
$X$ is nonsingular. 

The condition $\curl X = f X$ translates to
$d\iota_Xg = f\iota_X\mu$, or, $d\alpha = f\iota_X\mu$. 
Since $\mu$ is a volume form, its representation on each chart
is of the form $h e^1\wedge e^2\wedge e^3$, with $h$ a 
nonzero function. Hence, $\alpha$ is a contact form since
\begin{equation}
        \alpha\wedge d\alpha = \iota_Xg\wedge f \iota_X\mu
        = f h\norm{X}^2 e^1 \wedge e^2 \wedge e^3 \neq 0.
\end{equation} 
Finally, $X$ is Reeb-like with respect to $\alpha$ since
\begin{equation}
        \iota_X d\alpha = f \iota_X\iota_X\mu = 0.
\end{equation}

Conversely, assume further that $\alpha$ is a contact form for $M$ having
Reeb field $X$. Assume that $Y=hX$ for some $h>0$. Then,
on charts for $M$ choose a parallelization $\{e_i\}$ of $M$ 
such that $e_1=X$ and $e_2$ and $e_3$ are in $\xi$, 
the kernel of $\alpha$, and form a symplectic 
basis for this plane field (this is always possible since 
$\rest{d\alpha}{\xi}$ is a symplectic structure); hence, 
$d\alpha(e_2,e_3)=1$. Again let $\{e^i\}$ denote the dual 1-forms
to the $\{e_i\}$. On charts, choose the following metric adapted to the 
flow in the $\{e_i\}$-coordinates:
\begin{equation}
g=\left[\begin{array}{ccc}
        h\inv&0&0 \\ 0&1&0 \\ 0&0&1 
\end{array}\right]  .
\end{equation}
The transition maps respect this locally defined metric since
$e_1$ is globally defined and the symplectic basis $\{e_2, e_3\}$ 
differs from chart to chart by an element of $U(1)$, 
(the plane field is globally defined).
We remark that in this metric, $\norm{X}=1/\sqrt{h}$ 
and $\norm{Y}=\sqrt{h}$.

We claim that, under the metric $g$, $Y$ is volume
preserving and parallel to 
its curl. First, $\iota_Yg = e^1$ since for $i=1,2,3$,
\begin{equation}
        \iota_Yg(e_i) = g(Y,e_i) = g(he_1,e_i) = \delta_{i1} ,
\end{equation}
where $\delta_{ij}$ is the Kronecker delta. 
Next, note that $e^1=\alpha$, since they act on the basis
$\{e_i\}$ identically. Specifically, $\iota_{e_1}\alpha=\iota_X\alpha=1$
and $\iota_{e_2}\alpha=\iota_{e_3}\alpha=0$ since $e_2$ and $e_3$ 
were chosen to lie in $\xi$. 
We now have that $de^1=e^2\wedge e^3$ since for $i<j$,
\begin{equation}
        de^1(e_i,e_j) = d\alpha(e_i,e_j) = 
                \delta_{2i}\delta_{3j}-\delta_{2j}\delta_{3i},
\end{equation}
since we also chose the pair $(e_2,e_3)$ to form a symplectic 
basis for $\xi$. 

Let $\mu$ denote the volume form $h\inv e^1\wedge e^2\wedge e^3$. 
Then,
\begin{equation}
d\iota_Y\mu = d(\iota_{he_1}h\inv e^1\wedge e^2\wedge e^3) 
        = d(e^2\wedge e^3) = d^2(e^1) = 0 ;
\end{equation}
hence, $Y$ is volume preserving with respect to $\mu$. 
To show that $Y$ is Beltrami with respect to $g$ and $\mu$, 
it suffices to note that 
\begin{equation}
        d\iota_Yg = de^1 = e^2\wedge e^3 , 
\end{equation}
as well as 
\begin{equation}
        \iota_Y\mu = \iota_{he_1}h\inv e^1\wedge e^2\wedge e^3
                = e^2\wedge e^3.
\end{equation}
Hence, $\curl Y = Y$.
\qed

Note that $Y$ is divergence-free under the particular $g$-induced 
volume form if and only if the scaling function $h$ is an 
integral for the flow: \ie, $\L_Yh=0$.

\begin{rem}{\em
If $X$ is a Beltrami field with singularities, then we may
excise the singular points from the manifold and apply
Theorem~\ref{thm_BelReeb} to the nonsingular portion of the
flow. Then $X$ is still a Reeb flow for a contact form on 
the punctured manifold. 
}\end{rem}

\begin{cor}
Every Reeb-like field generates a steady-state solution to the 
Euler equations for a perfect incompressible fluid with respect 
to some metric.
\end{cor}

\begin{cor}\label{cor_Existence}
Every closed three-manifold has a rotational Beltrami flow on it 
for some metric. In fact, there are an infinite number, distinct up to 
homotopy through nonsingular vector fields. 
\end{cor}
\pf
The work of Martinet \cite{Ma} and Lutz \cite{Lu} shows that 
every closed three-manifold has a contact structure. 
To show that there 
are a countable collection of homotopically distinct Beltrami 
fields on any three-manifold, we appeal to the classification 
of plane fields on three-manifolds, elucidated recently by 
Gompf \cite{Go}. In the case where the rank of 
$H^2(M;\real)\neq 0$ (\ie, the cohomology is not pure 
torsional), there are at least a $\zed$'s worth of distinct 
Euler classes for plane fields on $M$; hence for 
(overtwisted) contact structures and the associated 
Beltrami fields. In the case where
$H^2(M;\zed)$ has only torsion elements, 
the {\em three-dimensional invariant} of Gompf \cite{Go} implies 
the existence of an infinite number of homotopy classes of plane 
fields: again, hence, of (overtwisted) contact structures and thus Beltrami
fields. 
\qed

Note that it is {\em a priori} unclear how one would construct 
a (nonzero) Beltrami flow on a nontrivial three-manifold, 
such as $(S^1\times S^2)\#(S^1\times S^2)$. Proving global
results about such flows would appear to be even more
unlikely. 

\section{Hofer's Theorem and the Hydrodynamical Seifert Conjecture}
\label{sec_WC}

A fundamental conjecture in global symplectic and contact
topology is the Weinstein Conjecture \cite{We}, which concerns 
certain embeddings of codimension-one submanifolds in
symplectic manifolds. We briefly recount this conjecture, following 
\cite{MS}. For excellent treatments, including recent progress, 
see the books \cite{HZ,MS}. 

Let $M^{2n-1}$ be a hypersurface in a symplectic manifold 
$(W^{2n},\omega)$.
Then $M$ is said to be of {\em contact type} if there exists 
a transverse expanding vector field $X$ on a neighborhood of 
$M$. That is, $X$ is transverse to $M$ and $\L_X\omega=\omega$.

\begin{conj}[Weinstein Conjecture, symplectic version \cite{We}]
\label{conj_WC}
Any hypersurface $M$ of contact type\footnote{Originally, the 
conjecture was stated under the additional hypothesis that 
$H^1(M;\rats)=0$. All known progress on this conjecture seems
to indicate that the cohomological condition is extranneous.}
has closed characteristics. That is, the foliation generated by 
the line field 
\[ 
\F = \left\{v\in T_zM :  \omega(v,T_zM)=0 \right\}, 
\]
has a closed leaf.
\end{conj}

Following \cite{MS}, it follows that any surface $M$ of contact type 
has a natural 1-form $\alpha$ associated to it via 
$\alpha=\iota_X\omega$, where $X$ is the expanding vector field. 
Since $d\alpha=\rest{\omega}{M}$, it follows that $\alpha$ is 
a contact form with Reeb field generating the foliation 
$\F$ above. Similarly, given a Reeb field on $M$, one can 
embed $M$ in a symplectic manifold as a hypersurface of 
contact type \cite{MS}. Hence, the Weinstein Conjecture translates to 
whether Reeb fields on $M$ have closed orbits. 
>From Theorem~\ref{thm_BelReeb}, we thus have:

\begin{conj}[Weinstein Conjecture, hydrodynamics version]
Every smooth\footnote{By ``smooth'' we mean $C^\infty$. However, the 
question is certainly interesting in other degrees of regularity.} 
rotational Beltrami flow on a Riemannian three-manifold has a 
closed orbit.
\end{conj}

The simplest example of an Euler flow without closed flowlines
is linear irrational flow on the Euclidean three-torus $T^3$. However, 
this flow has everywhere vanishing curl. A simple application 
of Stokes' Theorem implies that this flow (or any flow transverse to
a closed surface) cannot be rotational Beltrami under any metric. 

Great progress on the Weinstein Conjecture has been 
made in the past decade (see \cite{HZ} for a historical account). 
For example, Viterbo \cite{Vi} showed that 
the Weinstein Conjecture is true on any three-manifold of
contact type in $\real^4$ equipped with the standard symplectic form. 
More important for our applications, though, is the following 
seminal result of Hofer:

\begin{thm}[Hofer \cite{Ho}]\label{thm_Hofer}
If $\xi = \ker\alpha$ is an overtwisted contact structure on a compact 
orientable three-manifold, then the associated Reeb field
has a closed orbit of finite order in $\pi_1(M)$. 
\end{thm}

The proof of Theorem~\ref{thm_Hofer} is highly nontrivial, 
relying on the techniques of pseudo-holomorphic curves 
pioneered by Gromov \cite{Gr}. Such techniques yield, among other 
things, the following partial resolution to Conjecture~\ref{conj_WC}.

\begin{thm}[Hofer \cite{Ho}]\label{thm_HoferWC}
The Weinstein Conjecture is true for a three-manifold $M$
in the following cases: (1)  $M = S^3$ (or is covered by $S^3$); 
(2) $\pi_2(M) \neq 0$.
\end{thm}
In the case of $S^3$, the Reeb fields associated to overtwisted
structures are covered by Theorem~\ref{thm_Hofer}.
According to Eliashberg's classification of contact structures 
on $S^3$ \cite{El3}, any tight contact structure is contactomorphic
to that of Example~\ref{ex_TightS3}. The existence of a closed 
orbit then follows from techniques of Rabinowitz \cite{Ra}.

Theorem~\ref{thm_HoferWC} is the key step in the hydrodynamical
Seifert Conjecture:
\begin{thm}\label{thm_SeifertS3}
Any $C^\omega$ steady-state Euler flow on $S^3$ has a closed 
flowline.
\end{thm}
\pf
Given $X$ a nonsingular Euler field on $S^3$ under some metric $g$, 
Theorem~\ref{thm_Arnold} presents two possible cases. 

{\bf Case I:}
If $X$ is not everywhere colinear with $\curl X$, then the 
reduced pressure function $p$ is an integral for the flow. Hence, 
$S^3$ is decomposed by a compact subset $\Sigma$ into a finite 
collection of cells diffeomorphic to $T^2\times\real$, where the
each $T^2\times\{x\}$ is invariant under the flow. The slope of the 
vector field on each cell may be constant and irrational; hence, we must
consider the singular set, $\Sigma$, which is the inverse image of 
the (finite number of) critical values of the function $p:S^3\ra\real$ 
from Equation~\ref{eq_Euler}.

The key tool in this case comes from singularity theory: the inverse image 
of a critical value of $p$ is a (Whitney) stratified set, 
since it is a real-analytic variety \cite{GM}. That is, although the subset
is not necessarily a manifold, it can be decomposed into manifolds
of varying dimension ($\leq 2$ in our case) 
glued together in a sufficiently regular manner:
see \cite{GM} for proper definitions and theory. 

To determine the structure of $\Sigma$, 
note that it is locally a product in the flow direction. More 
specifically, let $x\in\Sigma$. Since $X$ is nonsingular, the 
Flowbox Theorem (see \eg, \cite{Ro}) states that there is a neighborhood 
$U\cong D^2\times\real$ of $x$ in $S^3$ such that $\rest{X}{U}$ 
points entirely in the $\real$-component. Denote by $D_0$
the disc $D^2\times\{0\}$ containing $\{x\}$. Since $S^3\setminus\Sigma$ 
is filled with invariant tori, $\Sigma$ is invariant and the Flowbox
Theorem implies that $D_0$ is transverse to $\Sigma$ and
$U\cap\Sigma\cong(D_0\cap\Sigma)\times\real$. Thus, the transverse 
structure of $\Sigma$ is invariant along flowlines of $X$.

By analyticity, $\Sigma$ is at most two-dimensional. This, along 
with the Whitney condition (a) for stratified sets \cite{GM}, 
implies that the intersection $D_0\cap\Sigma$ is homeomorphic 
to a radial $k$-pronged tree centred at $x$. 
If the number $k$ is zero, then $\Sigma$ is a 1-manifold, which by the 
compactness of $\Sigma$ and the flow-invariance of the transverse
structure, implies that $X$ has a compact invariant 1-manifold. 
If $k=1$, then $\Sigma$ has an invariant 
1-dimensional boundary component, which again by compactness 
and flow-invariance implies a closed orbit.
If $k>2$, then, by the Uniqueness Theorem for 
ODEs, $x$ lies within an invariant 1-stratum, which as before 
must continue to a compact invariant 1-manifold via flow-invariance 
and finiteness of the stratification. 

Thus, if $X$ has no closed orbits on $\Sigma$, then for 
every $x\in\Sigma$, the transverse prong number $k$ is precisely 
two, and $\Sigma$ is locally homeomorphic to $\real^2$. 
Hence it is an invariant closed 2-manifold with nonsingular vector
field: $\Sigma=\coprod T^2$. Thus, $S^3$ is 
decomposed into a finite collection of $T^2\times\real$ cells 
glued together (pairwise) along sets diffeomorphic to $T^2$, and 
we have expressed the three-sphere as a $T^2$-bundle 
over $S^1$, a contradiction. 

{\bf Case II:}
If $X$ is everywhere colinear with $\curl X$, then 
$\curl X=f X$ and $X$ is a Beltrami field. 
In the case where $f$ is nonzero, $X$ is rotational and  
the analysis of Case I is completely useless: unlike all other
cases, we have a priori no information about the flow. However, 
by Theorem~\ref{thm_BelReeb}, this vector
field is, up to a rescaling, a Reeb field. Hence, the associated flows
are related by a time reparametrization, which neither 
creates nor destroys the closed orbit whose existence follows 
from Theorem~\ref{thm_HoferWC} of Hofer.  

If $f$ has zeros, then either $f$ varies, or it is constantly zero.
If $f$ is not constant, then $f$ is a nontrivial integral 
for $X$ as follows. First, $\curl X$ is $\mu$-preserving since
\begin{equation}
        d\iota_{\curl X}\mu = d(d\iota_X g) = 0.
\end{equation}
However, the left hand term simplifies to 
\begin{equation}
        d\iota_{\curl X}\mu = d(f\iota_X\mu) = df\wedge\iota_X\mu ,
\end{equation}
since $X$ is also $\mu$-preserving. Finally, since we are 
on a 3-manifold, the 4-form $df\wedge\mu$ vanishes, implying
\begin{equation}
        0 = \iota_X(df\wedge\mu) = (\iota_X df)\mu - df\wedge\iota_X\mu 
                = (\iota_X df)\mu. 
\end{equation}
Thus, $\iota_X df=0$ and $f$ is an integral. 
The existence of a closed orbit then 
follows from the singularity-theory arguments of Case I. 

Finally, if $\curl X = 0$ everywhere, then the 1-form 
$\alpha=\iota_X g$ is closed and nondegenerate since 
$d\iota_Xg = f\iota_X\mu = 0$. By the Frobenius condition, 
$\xi=\ker\alpha$ defines a smooth codimension-one foliation
of $S^3$, to which $X$ is transverse. However, by the Novikov
Theorem \cite{No}, $\xi$ has a Reeb component: in 
particular, there is a closed leaf homeomorphic to $T^2$. 
This torus separates in $S^3$ and hence cannot be transverse to 
the volume-preserving field $X$, contradicting the assumption 
that $X$ is nonsingular.  
\qed

\begin{rem}{\em
For a typical closed orientable three-manifold $M$, the Weinstein 
Conjecture is still open. However, we may conclude that 
``most'' Beltrami fields (in the homotopic sense) possess closed 
flowlines: by Corollary~\ref{cor_Existence},
there are an infinite number of homotopically distinct 
rotational Beltrami fields on $M$. By a theorem of Kronheimer and
Mrowka \cite{KM}, there are on $M$ only a finite number of 
homotopy classes of plane fields which contain a tight contact
representative. Hence, by Theorems~\ref{thm_BelReeb} 
and \ref{thm_Hofer}, all but finitely many of the rotational Beltrami 
fields on $M$ (up to homotopy) have closed flowlines.
}\end{rem}

The important feature of Hofer's Theorem is the relationship
between the ``twistedness'' of a contact structure and 
the dynamics transverse to it. In order to apply 
contact-topological methods to manifolds which are of greater interest to 
fluid dynamicists ($T^3$ or $\real^3$) we consider 
algebraic obstructions to tightness.

\section{Beltrami flows on $T^3$} \label{sec_ABC}

The most important three-manifolds from the point of view
of hydrodynamics are the solid torus $D^2\times S^1$ and 
the three-torus $T^3=S^1\times S^1\times S^1$. Some authors 
have suggested using Beltrami ``tubes'' to model certain 
domains of turbulent regions in flows, following obervations 
that suggest the vorticity tends to align with velocity in 
certain domains \cite{Mo}. We do not 
discuss the applications of contact geometry to this 
case, since little is known of the classification of tight 
contact structures on $D^2\times S^1$. 

Contact structures on the three-torus, however, are 
more completely understood. In addition, the most important 
examples of interesting Beltrami fields live on $T^3$. 
A fundamental observation of Arnold's is the existence of
Beltrami fields on the three-torus $T^3$ which are nonintegrable: 
these so-called {\em ABC flows} (and their generalizations) 
have been the source of a great deal of inquiry:
\begin{equation}\label{eq_ABC}
\begin{array}{l}
\dot{x} = A\sin z + C\cos y \\
\dot{y} = B\sin x + A\cos z \\
\dot{z} = C\sin y + B\cos x 
\end{array}
,\end{equation}
for some $A, B, C \geq 0$. By symmetry in the variables 
and parameters, we may assume without loss of generality that 
$1=A\geq B\geq C\geq 0$. Under this convention,
the vector field is nonsingular if and only if \cite{D+}
\begin{equation}\label{eq_Nonsingular}
        B^2 + C^2 \leq 1.
\end{equation}

Though the list of publications concerning ABC flows is 
extensive, there is very little known about the global features 
of these flows, apart from cases where one of the constants
(say $C$) is zero or a perturbation thereof. Melnikov 
methods 
have been used repeatedly to show complex
behavior in these near-integrable cases. Of the generalizations
of ABC flows to other eigenfields of the curl operator, even 
less is known. There are seemingly no global results on 
general Beltrami fields.
However, thanks to the following classification theorem of 
Giroux, we may apply our contact-topological
methods to the most general case of Beltrami flows on $T^3$:

\begin{thm}[Giroux \cite{Gi}] \label{thm_Giroux}
Any tight contact structure on $T^3$ is (up to a choice of 
fundamental class in $H_2(T^3)$) isotopic through contact structures 
to the kernel of $\alpha_n$ for some integer $n>0$, where
\begin{equation}\label{eq_Giroux}
        \alpha_n = \sin(nz)dx + \cos(nz)dy .
\end{equation}
\end{thm}

\begin{dfn}\label{def_Homotopy} {\em
Let $X$ be a nonsingular vector field on $M$. 
Since all three-manifolds are parallelizable, $X$ 
gives a map from $M$ to $\real^3$.
Under the normalization $\real^3\setminus\{0\}\ra S^2$, $X$ 
induces a map $M\ra S^2$. The {\em homotopy
class} of $X$ is defined to be the homotopy class 
of the map $M\ra S^2$.
}\end{dfn}

In other words, a nonsingular vector field
is homotopically trivial if it can be deformed through nonsingular
vector fields in such a way that all the vectors ``point in the 
same direction.''

\begin{thm}\label{thm_T3}
Every homotopically nontrivial rotational Beltrami field on $T^3$ 
has contractible closed orbits.
\end{thm}
\pf
Let $X$ denote a homotopically nontrivial field on a 
Riemannian $T^3$ satifying $\curl X=fX$ with $f>0$. 
Then by Theorem~\ref{thm_BelReeb}, there is a natural contact
structure $\xi$ transverse to $X$ and uniquely defined up to 
homotopy. If $X$ is homotopically nontrivial, then so is 
$\xi$ as a plane field, since the homotopy class of an oriented 
plane field is defined as the homotopy class of the vector field 
transverse to it. 

By Theorem~\ref{thm_Giroux}, any tight contact structure
on $T^3$ is isotopic to the kernel of the 1-form
$\sin(nz)dx + \cos(nz)dy$. The Reeb field $X_n$ 
associated to this contact form $\alpha_n$ is
\begin{equation}\label{eq_T3Reeb}
        X_n = \sin(nz)\frac{\del}{\del x} + \cos(nz)\frac{\del}{\del y}.
\end{equation}
As the image of the induced map $T^3\ra S^2\subset\real^3$ lies on 
the equator $z=0$, it follows that every tight contact structure
on $T^3$ is homotopically trivial.
Hence, $X$ is a Reeb-like field for an overtwisted contact structure,
which, by Theorem~\ref{thm_Hofer}, must have 
a closed orbit which is torsional in $\pi(T^3)\cong\zed^3$ --- \ie,
it is contractible.
\qed

\begin{rem}{\em
The existence of closed orbits which are {\em contractible} is 
of particular importance. One typically works on $T^3$ in order
to model spatially periodic flows on $\real^3$. The existence 
of a contractible closed orbit on $T^3$ implies that when lifted
to the universal cover $\real^3$, the orbit remains closed. We 
note that upon numerically integrating examples of Beltrami 
fields on $T^3$ (\eg, the ABC equations), 
the integrable regions are certainly homotopically nontrivial 
in $T^3$, whereas the closed orbits in the nonintegrable regions
are completely obscured.
}\end{rem}

\begin{rem}{\em
The determination of the homotopy type of a nonsingular 
vector field on a three-manifold is a problem of algebraic topology,
thanks to recent results of Gompf \cite{Go}. Gompf assigns to 
a vector field a pair of integer-valued invariants: a two-dimensional
refinement of the Euler class, and a more subtle 
three-dimensional invariant, which is derived from what type 
of four-manifold $M$ bounds. Together, these invariants, which 
can be computed in many cases, completely classifiy
the homotopy type. 
}\end{rem}

\begin{rem}{\em
It is by no means the case that Theorem~\ref{thm_T3} holds 
for vector fields in general. 
As mentioned earlier, the results of Kuperberg \cite{Ku1} allow 
one to insert ``plugs'' to break isolated closed orbits. Since the 
Kuperberg plugs do not change the homotopy type of the 
vector field, there are smooth nonsingular vector fields on $T^3$ 
in {\em every} homotopy class which have no closed orbits. 
}\end{rem}

We may extend Theorem~\ref{thm_T3} to the analogue of the
Seifert Conjecture for Euler flows on $T^3$.
\begin{thm}\label{SeifertT3}
Any homotopically nontrivial $C^\omega$ Euler flow on $T^3$ has a 
closed orbit. 
\end{thm}
\pf
By Theorems~\ref{thm_T3}, \ref{thm_Arnold}, and the proof
of Theorem~\ref{thm_SeifertS3}, the only 
remaining cases are when one has an integrable Eulerian field $X$
on $T^3$, or when $\curl X = 0$. 

Consider first the integrable case.
As in the proof of Theorem~\ref{thm_SeifertS3}, we
know that either there exists a closed invariant 1-manifold, or we have
constructed $T^3$ as a union of cells $T^2\times\real$ glued
together along invariant 2-tori. In the latter case, if there are 
no closed orbits, then the rotation number for the flow on each 
$T^2$ is irrational and thus constant. Hence, we have a 
homotopically trivial flow, since each $T^2\times\{x\}$ is 
invariant and has linear flow on it (up to conjugacy), so that the 
image of the Gauss map is contractible in $S^2$.

In the remaining case where $\curl X=0$, we again know that 
the 1-form $\alpha=\iota_Xg$ defines a smooth codimension 
one foliation which is taut since $X$ is a transverse volume-preserving
vector field. However, smooth taut foliations can be perturbed
to a tight contact structure \cite{ET}. Thus, by 
Theorem~\ref{thm_Giroux}, the foliation (and the corresponding
normal vector field) is homotopically trivial.
\qed

\begin{rem}{\em
One may also utilize the correspondence of Theorem~\ref{thm_BelReeb}
to construct interesting examples of Reeb fields. For example, 
the analysis of Dombre {\em et al.} \cite{D+} indicates the
presence of integrable ``tubes'' within the ABC flows. One 
may cut out these tubes and replace them in a topologically 
different manner in such a way as to produce a new contact 
structure on a different three-manifold. This form of Dehn 
surgery which is adapted to Reeb fields can be done so that 
the new Reeb field agrees with the ABC field outside of the 
integrable tubes \cite{EG}: hence the flow is 
nonintegrable and exhibits Lagrangian turbulence. 
Since the integrable tubes present in the ABC flows form a 
complete set of generators for the first homology group 
$H_1(T^3;\zed)$, we can break all the first homology and 
obtain an ABC-type flow on the three-sphere $S^3$. 
}\end{rem}

We close with the specific case of the ABC flows. Unfortunately, 
we cannot apply Theorem~\ref{thm_T3} to these equations.
\begin{prop}
Every nonsingular ABC field is transverse to a tight contact 
structure.
\end{prop}
\pf
By normailizing the coefficient $A$ to 1 and using 
Equation~\ref{eq_Nonsingular}, we have that the parameter space
$\{(B,C): 0\leq B^2+C^2\leq 1\}$ is path-connected; hence, 
if we show that
some ABC field satisfies the proposition, then every other ABC 
field is homotopic to this through nonsingular Beltrami fields, and 
so the transverse contact structures are isotopic through 
contact structures. In the particular case where $B=C=0$, we
have the equations
\begin{equation}\begin{array}{l}
        \dot{x} = A\sin z \\
        \dot{y} = A\cos z
\end{array} ,
\end{equation}
which is the Reeb field for the (tight) contact form 
$\alpha = A\sin z dx + A\cos z dy$ (\cf Equation~\ref{eq_T3Reeb}). 
\qed

\begin{rem}{\em
An important feature of any Beltrami field $X$ is the fact that it 
extremizes the energy functional 
\begin{equation}\label{eq_Energy}
        E(\tilde{X}) = \frac{1}{2}\int_M\norm{\tilde{X}}^2 d\mu
\end{equation}
among the class of all vector fields $\tilde{X}$ obtained from $X$ by 
$\mu$-preserving diffeomorphisms of $M$ \cite{Ar2}. It 
would be interesting to investigate the relationship between 
smooth fields which {\em minimize} the energy functional and the
associated transverse contact structures. It follows from 
remarks in Arnold \cite{Ar2} that the Reeb field 
associated to the unique tight contact structure on $S^3$, as 
well as the ABC flows each minimize energy. In fact, 
every known example of a smooth energy-minimizing 
field is the Reeb field for a tight contact structure. It would 
be interesting to see whether one could always reduce
the energy of a Reeb field associated to an overtwisted 
contact structure: if so, this may be a significant application of 
the tight / overtwisted dichotomy in contact geometry to a problem
in hydrodynamics.
}\end{rem}

\vspace{0.1in}

\Large
\bf{Acknowledgements}
\normalsize

The authors acknowledge the support of the University of Texas 
at Austin. This work is also supported by a National Science
Foundation Postdoctoral Fellowship [RG], number DMS-9508846. 
We were greatly aided through conversations with 
Emmanuel Giroux, Rafael de la Llave, Robert MacPherson,
Igor Mezi\'c, Karen Uhlenbeck, Misha Vishik, and Bob Williams.  
Rafael de la Llave and Margaret Symington carefully read
a manuscript and provided helpful clarifications.


\end{document}